\documentclass[10pt,twoside]{article}
\usepackage{graphicx}
\usepackage{amsmath}
\usepackage{amssymb}
\usepackage{Latex-document}
\usepackage{titlesec}

\newcommand{\hk}{\hat{k}}

\newcommand{\ZZo}{\Z^2 /\{0\}}

\newcommand{\Z}{{\mathbb Z}}
\newcommand{\Zo}{\Z /\{0\}}

\newcommand{\e}{\epsilon}
\newcommand{\vth}{\vartheta}

\newcommand{\C}{{\mathcal C}}
\renewcommand{\k}{\kappa}
\newcommand{\ga}{\gamma}

\newcommand{\Dl}{\Delta}
\renewcommand{\th}{\theta}

\newcommand{\ra}{\rightarrow}
\newcommand{\al}{\alpha}
\newcommand{\be}{\beta}

\newcommand{\pa}{\partial}

\newcommand{\la}{\lambda}

\newcommand{\om}{\omega}
\newcommand{\Om}{\Omega}
\newcommand{\na}{\nabla}

\newcommand{\tq}{\tilde{q}}

\newcommand{\cO}{{\mathcal O}}
\newcommand{\non}{\nonumber}

\newcommand{\tom}{\tilde{\omega}}
\newcommand{\R}{\mathbb{R}}

\renewcommand\thesection{\arabic{section}}
\titleformat{\section}{\normalfont\Large\bfseries}{\thesection.}{0.5em}{}

\title{\bf Chaos in Partial Differential Equations, Navier-Stokes Equations and Turbulence}

\author{Y. Charles Li \thanks{Department of Mathematics,
University of Missouri, Columbia, MO65211, USA. E-mail: cli@math.missouri.edu}}
\date{}

\begin{document}
\maketitle

\thispagestyle{first} \setcounter{page}{1}

\begin{abstract}\vskip 3mm\footnotesize
\noindent I will briefly survey the most important results obtained so far on 
chaos in partial differential equations. I will also survey progresses and 
make some comments on Navier-Stokes equations and turbulence.

\vskip 4.5mm


\end{abstract}

\vskip 12mm

\section{Introduction} 

In my view, the two most important problems in the area of chaos (turbulence) are
\begin{enumerate}
\item An effective description of chaotic (turbulent) solutions beyond the Reynolds average,
so that important characteristics of the chaotic (turbulent) solutions can be revealed 
accurately.
\item Proofs of the existence of chaos (turbulence), so that the parameter regime and 
certain characteristics of chaos (turbulence) can be discovered. Control of chaos (turbulence)
is a sub-problem here.
\end{enumerate}
So far, there is no substantial progress on problem 1 \cite{Li07a}. There are 
important progresses on problem 2 \cite{Li07b} which will be addressed later. 
A notable mathematical problem on Navier-Stokes equations is the global regularity problem
which was chosen as one of the seven Clay problems. 

Extensive studies have been conducted on chaos in ordinary differential equations
(ODE) \cite{ASY97}, and intensive studies are focused upon Lorenz equations 
\cite{Spa82}. Existence of chaos in Lorenz equations has not been proved by hand. 
Recently, there have been a new trend of computer proofs \cite{Tuc02}. In general,
for near integrable ODE with transversal homoclinics, existence of chaos can 
often be proved by hand. The same is true for PDE. 

Chaos in partial differential equations (PDE) has long been an open area. During 
the last decade, a standard program was established for proving the existence of chaos
in near integrable PDE \cite{Li04}. Around transversal homoclinics, existence of 
chaos can be proved by hand \cite{Li03a} \cite{Li03b} \cite{Li04c}, while around 
non-transversal homoclinics, existence of chaos can be proved by hand up to 
nasty generic conditions \cite{Li04a} \cite{Li04b}. There are also attempts 
of computer proofs on the existence of chaos in PDE. 

In contrast to ODE, PDE has a lot of novelties. For instance, boundary condition 
plays a major role on dynamics. So far, all the chaos in PDE is proved under 
periodic boundary condition. In fact, odd or even constraint besides the periodic 
boundary condition is crucial. Often one can prove the existence of chaos under 
odd or even constraint, but not in general. Periodic boundary condition permits a 
representation of solutions by Fourier series. The sequence of Fourier modes in the 
Fourier series leads to a natural generalization of ODE to infinite dimensional 
systems. The Fourier series also makes hand calculation possible. Under other 
boundary conditions, hand calculation is often impossible. More importantly, the 
dynamics is dramatically different. For example, for the whole space problem 
(i.e. under decaying boundary conditions), existence of chaos has not been proved.
Of course, the PDE type is fundamental to dynamics. So far, the chaos proved 
is for mixed hyperbolic and parabolic semilinear systems which are near 
hyperbolic semilinear integrable systems. For semilinear systems, existence of
invariant manifolds is easy to obtain. For hyperbolic quasilinear systems, existence 
of invariant manifolds is a big open problem. Here one of the interesting examples 
will be Euler equations. Even for the simple integrable 
derivative nonlinear Schr\"odinger equation (DNLS)
\begin{equation}
iq_t = q_{xx} -i \left ( |q|^2 q \right )_x \ ,
\label{DNLS}
\end{equation}
existence of invariant manifolds is open while its invariant subspaces are perfectly 
normal \cite{Li07c}. This difficulty is the major obstacle toward proving the 
existence of chaos when the DNLS is under perturbations. 

Chaos in PDE is a rich area, full of potential, with a crown goal of solving the 
problem of turbulence. Turbulence is governed by Navier-Stokes equations (NS). Turbulence
often happens when the Reynolds number is large. Formally setting the Reynolds number
to infinity in the Navier-Stokes equations, one gets the Euler equations. Euler 
equations are a lot like hyperbolic quasilinear integrable systems \cite{Li01} 
\cite{LY03} \cite{Li03}. By posing 
periodic boundary condition to Navier-Stokes equations, one can study turbulence 
in the near Euler equations regime as a canonical problem. So far, there is no 
success in proving the existence of chaos in Navier-Stokes equations. 

\section{Chaos in Partial Differential Equations}  

The standard program that we developed for proving the existence of chaos in partial 
differential equations \cite{Li04} involves multi-disciplinary subjects including
integrable theory, dynamical system, partial differential equation, and functional analysis.
The specific machineries are: (1). Darboux transformations, (2). Isospectral theory, 
(3). Persistence of invariant manifolds and Fenichel fibers, (4). Melnikov analysis and shooting 
technique, (5). Shadowing technique and symbolic dynamics, (6). Specific Smale horseshoe 
construction.

\subsection{Transversal Homoclinics and Heteroclinics}

In a non-autonomous system (e.g. a periodic system), transversality can often be signified 
by its Poincar\'e period map. A typical example in this case for a success of proving 
the existence of chaos is the sine-Gordon system \cite{Li04c} \cite{Li06}:
\begin{equation}
u_{tt}=c^2 u_{xx} + \sin u +\e [ -a u + f \sin^3 u ]
\label{PSG}
\end{equation}
which is subject to periodic boundary condition and even (or odd) constraint
\begin{equation}
u(t, x+2\pi ) = u(t, x), \ u(t, -x) = u(t, x), \text{ or } u(t, -x) = - u(t, x),
\label{BC}
\end{equation}
where $u$ is a real-valued function of two real variables $t$ 
and $x$, $c$ and $a >0$ are parameters, $\frac{1}{2} < c < 1$, $\e$ is a small 
perturbation parameter, $\e \geq 0$, and $f = \cos t$ for example.

{\bf Theorem 1.} \cite{Li04c} \cite{Li06} \it For ($\e , a$) $\in$ ($0, \e_0$] $\times$ ($a_1,a_2$), 
where $\e_0 >0$ and $0<a_1 <a_2$, under the odd constraint (\ref{BC}), there is a homoclinic orbit; 
under the even constraint (\ref{BC}), there is a heteroclinic cycle. In a neighborhood of the              
homoclinic orbit or heteroclinic cycle, there is chaos. That is, there is a Cantor set $\Xi$ of 
points, which is invariant under 
an iteration of the Poincar\'e period map $F^K$ for some $K$. The action of $F^K$ on $\Xi$ is 
topologically conjugate to the action of the Bernoulli shift on two symbols $0$ and $1$. \rm

The proof on the existence of a homoclinic orbit or a heteroclinic cycle involves 
Darboux transformations, isospectral theory, persistence of invariant manifolds and Fenichel fibers,
and Melnikov analysis and shooting technique. The proof on the existence of chaos involves 
shadowing technique and symbolic dynamics.

The concept of a homoclinic orbit or a heteroclinic cycle can be generalized to a homoclinic tube or 
a heteroclinically tubular cycle \cite{Li04c} \cite{Li06} by replacing orbits with invariant tubes. 
It turns out that the above theorem is still true for tubes. On the other hand, dynamics inside the 
tubes can be chaotic too. So we have a smaller scale chaos embedded inside a larger scale chaos. In 
principle, this process can continue to smaller and smaller scale chaos. Thus we can have a 
``chaos cascade'' \cite{Li04c} \cite{Li06}. For example, let $f = f(t, \th )$ in (\ref{PSG}) where
$f(t, \th )$ is periodic in $t$ and $\th \in {\mathbb T}^N$. For 
example, $\th = (\th_1, \th_2, \th_3)$, $\th_n =\om_n t + \th_n^0
+\e^\mu \vth_n$, $\mu > 1$, and $\vth_n$'s are given by the ABC flow \cite{DFGHMS86},
\begin{eqnarray}
\dot{\vth}_1 &=& A \sin \vth_3 + C \cos \vth_2 \ , \non \\
\dot{\vth}_2 &=& B \sin \vth_1 + A \cos \vth_3 \ , \label{ABC} \\
\dot{\vth}_3 &=& C \sin \vth_2 + B \cos \vth_1 \ , \non 
\end{eqnarray}
which is chaotic for certain values of the real parameters $A$, $B$, and $C$.
Let $\th_4 = \om_4 t + \th_4^0$, $f$ can be, for example, of the form
\begin{equation}
f(t, \th )= \al +\sum_{n=1}^4 \be_n \cos \th_n \ ,
\label{force}
\end{equation}
where $\om_n$'s ($1\leq n \leq 4$) form a quasiperiodic basis, 
$\th_n^0$'s ($1\leq n \leq 4$) are constant phases, and $\al$ and $\be_n$ 
are real constants. 

{\bf Theorem 2.} \cite{Li04c} \cite{Li06} \it For ($\e , a$) $\in$ ($0, \e_0$] $\times$ ($a_1,a_2$), 
where $\e_0 >0$ and $0<a_1 <a_2$, under the odd constraint (\ref{BC}), there is a homoclinic tube; 
under the even constraint (\ref{BC}), there is a heteroclinically tubular cycle. In a neighborhood of 
the homoclinic tube or heteroclinically tubular cycle, there is tubular chaos. That is, there is a 
Cantor set $\Xi$ of tori, which is invariant under 
an iteration of the Poincar\'e period map $F^K$ for some $K$. The action of $F^K$ on $\Xi$ is 
topologically conjugate to the action of the Bernoulli shift on two symbols $0$ and $1$. \rm

In an autonomous system, a transversal homoclinic orbit or a heteroclinic cycle is asymptotic to a 
limit cycle. A typical example in this case for a success of proving the existence of chaos is the
Ginzburg-Landau equation in the near nonlinear Schr\"odinger regime:
\begin{equation}
i q_t = q_{xx} + 2 |q|^2 q +i \e \bigg [ \left (\frac{9}{16}-|q|^2 \right )q +\mu 
|\hat{\pa}_x q|^2 \bar{q} \bigg ]\ , \label{derNLS}
\end{equation}
where $q$ is a complex-valued function of two real variables $t$ and $x$,
$\e \geq 0$ is the perturbation parameter, $\mu$ is a real constant, and
$\hat{\pa}_x $ is a bounded Fourier multiplier,
\[
\hat{\pa}_x q = -\sum_{k=1}^K k \tq_k \sin kx\ , \quad 
\mbox{for} \ q = \sum_{k=0}^\infty \tq_k \cos kx\ ,
\]
for some fixed large $K$. Periodic boundary condition 
and even constraint are imposed,
\[
q(t,x+2\pi ) = q(t,x)\ , \ \ q(t,-x)=q(t,x) \ . 
\] 

{\bf Theorem 3.} \cite{Li03a} \it There exists a $\e_0 > 0$, such that 
for any $\e \in (0, \e_0)$, and $|\mu | > 5.8$,
there exist two transversal homoclinic orbits asymptotic to 
the limit cycle $q_c = \frac{3}{4} \exp \{ -i [ \frac{9}{8} t + \ga ]\}$. 
In a neighborhood of the homoclinic orbits there is chaos, in the sense of Theorem 1 
with the Poincar\'e period map replaced by a Poincar\'e return map. \rm

The proof on the existence of the two transversal homoclinic orbits uses the same tools as above. 
The proof on the existence of chaos uses a shadowing lemma developed in \cite{Li03a} for 
infinite dimensional autonomous systems --- a long time open problem. 

\subsection{Non-Transversal Homoclinics and Heteroclinics}

In an autonomous system, homoclinics or heteroclinics asymptotic to a saddle is 
usually non-transversal. Some of such homoclinics or heteroclinics can still induce chaos,
e.g. Silnikov homoclinics. Nevertheless, a proof without any generic assumption on such 
chaos is still elusive. A typical example is the Ginzburg-Landau equation in the near 
nonlinear Schr\"odinger regime \cite{Li04a}:
\begin{equation}
iq_t = q_{xx} +2 [ |q|^2 - \om^2] q +i \e [q_{xx} - \al q +\be ] \ ,
\label{pnls}
\end{equation}
where $q = q(t,x)$ is a complex-valued function of the two real 
variables $t$ and $x$, $\om$ $\in$ ($1/2 , 1$), $\al >0$ and $\be >0$ are constants, 
and $\e \geq 0$ is the perturbation parameter. Periodic boundary condition 
and even constraint are imposed,
\[
q(t,x + 2 \pi) = q(t,x)\ , \ \ q(t,-x) = q(t,x)\ .
\]

\begin{figure}[ht]
\includegraphics[width=4.0in,height=3.0in]{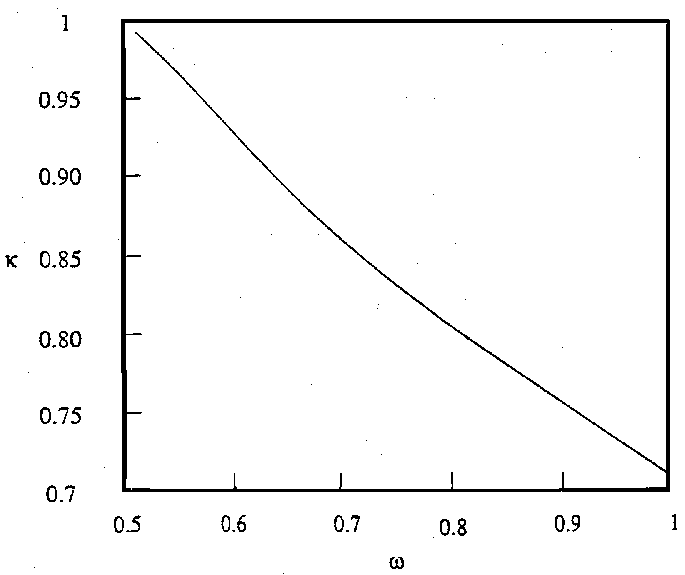}
\caption{The graph of $\k =\k (\om )$.}
\label{kap}
\end{figure}

{\bf Theorem 4.} \cite{Li04a} \cite{Li04b} \it There exists a $\e_0 > 0$, such that 
for any $\e \in (0, \e_0)$, there exists 
a codimension 1 surface in the space of $(\alpha,\beta, \om) \in 
\mathbb{R}^+\times  \mathbb{R}^+\times \mathbb{R}^+$ where 
$\om \in (\frac{1}{2}, 1)/S$, $S$ is a finite subset, and 
$\al \om < \be$. For any $(\alpha ,\beta, \omega)$ on the codimension-one
surface, the Ginzburg-Landau equation (\ref{pnls}) has a pair of Silnikov homoclinic orbits 
asymptotic to a saddle. The codimension 1 surface has the approximate expression 
$\al = 1/ \k (\om )$ where $\k (\om )$ is given in Figure \ref{kap}. Under certain 
generic assumptions, Smale horseshoes can be constructed in the neighborhood of the 
homoclinic orbits. That is, there exists chaos under generic assumptions. \rm

\subsection{Extra Comments}

Other success or partial success in proving the existence of chaos includes 
the discrete Ginzburg-Landau equation in the near discrete nonlinear Schr\"odinger regime
\cite{Li92} \cite{LM97} \cite{LW97} \cite{Li03t}, perturbed Davey-Stewartson equation 
\cite{Li00b} \cite{Li05b}, perturbed vector nonlinear Schr\"odinger equation \cite{Li02b}. 

Arnold diffusion can be proved for discrete nonlinear Schr\"odinger equation under certain 
Hamiltonian perturbation \cite{Li06z}.

\section{Navier-Stokes Equations and Turbulence}

The success in proving the existence of chaos in near integrable partial differential equations 
naturally leads to the question of proving the existence of chaos in Navier-Stokes equations (NS)
in the high Reynolds number (near Euler equations) regime. Like Lorenz system among systems of 
ODE, chaos in Navier-Stokes equations gathers the most interest. Unlike Lorenz system, 
high Reynolds number Navier-Stokes equations are a lot like near integrable systems, and 
Euler equations are a lot like integrable systems.

\subsection{Lax Pairs of Euler Equations}

The 2D Euler equation can be written in the vorticity form,
\begin{equation}
\pa_t \Om + \{ \Psi, \Om \} = 0 \ ,
\label{euler}
\end{equation}
where the bracket $\{\ ,\ \}$ is defined as
\[
\{ f, g\} = (\pa_x f) (\pa_y g) - (\pa_y f) (\pa_x g) \ ,
\]
$\Om$ is the vorticity, and $\Psi$ is the stream function given by,
\[
u=- \pa_y \Psi \ ,\ \ \ v=\pa_x \Psi \ ,
\]
and the relation between vorticity $\Om$ and stream 
function $\Psi$ is,
\[
\Om =\pa_x v - \pa_y u =\Dl \Psi \ .
\]

{\bf Theorem 5.} \cite{Li01} \it The Lax pair of the 2D Euler equation (\ref{euler}) is given as
\begin{equation}
\left \{ \begin{array}{l} 
L \varphi = \la \varphi \ ,
\\
\pa_t \varphi + A \varphi = 0 \ ,
\end{array} \right.
\label{laxpair}
\end{equation}
where
\[
L \varphi = \{ \Om, \varphi \}\ , \ \ \ A \varphi = \{ \Psi, \varphi \}\ ,
\]
and $\la$ is a complex constant, and $\varphi$ is a complex-valued function.\rm 

Starting from a Lax pair, one can try to develop a Darboux transformation which is the key in the 
success in proving the existence of chaos in near integrable PDE. Some progress is made in 
establishing the Darboux transformation for 2D Euler equation.
Consider the Lax pair (\ref{laxpair}) at $\la =0$, i.e.
\begin{eqnarray}
& & \{ \Om, p \} = 0 \ , \label{d1} \\
& & \pa_t p + \{ \Psi, p \} = 0 \ , \label{d2} 
\end{eqnarray}
where we replaced the notation $\varphi$ by $p$.

{\bf Theorem 6.} \cite{LY03} \it Let $f = f(t,x,y)$ be any fixed solution to the system 
(\ref{d1}, \ref{d2}), we define the Gauge transform $G_f$:
\begin{equation}
\tilde{p} = G_f p = \frac {1}{\Om_x}[p_x -(\pa_x \ln f)p]\ ,
\label{gauge}
\end{equation}
and the transforms of the potentials $\Om$ and $\Psi$:
\begin{equation}
\tilde{\Psi} = \Psi + F\ , \ \ \ \tilde{\Om} = \Om + \Dl F \ ,
\label{ptl}
\end{equation}
where $F$ is subject to the constraints
\begin{equation}
\{ \Om, \Dl F \} = 0 \ , \ \ \ \{ \Dl F, F \} = 0\ .
\label{constraint}
\end{equation}
Then $\tilde{p}$ solves the system (\ref{d1}, \ref{d2}) at 
$(\tilde{\Om}, \tilde{\Psi})$. Thus (\ref{gauge}) and 
(\ref{ptl}) form the Darboux transformation for the 2D 
Euler equation (\ref{euler}) and its Lax pair (\ref{d1}, \ref{d2}). \rm

For KdV equation and many other soliton equations, the 
Gauge transform is of the form \cite{MS91},
\[
\tilde{p} =  p_x -(\pa_x \ln f)p \ .
\]
In general, Gauge transform does not involve potentials.
For 2D Euler equation, a potential factor $\frac {1}{\Om_x}$
is needed. From (\ref{d1}), one has
\[
\frac{p_x}{\Om_x} = \frac{p_y}{\Om_y} \ .
\]
The Gauge transform (\ref{gauge}) can be rewritten as
\[
\tilde{p} = \frac{p_x}{\Om_x} - \frac{f_x}{\Om_x} \frac{p}{f}
=\frac{p_y}{\Om_y} - \frac{f_y}{\Om_y} \frac{p}{f}\ .
\]
The Lax pair (\ref{d1}, \ref{d2}) has a symmetry, i.e. it is 
invariant under the transform $(t,x,y) \ra (-t,y,x)$. The form 
of the Gauge transform (\ref{gauge}) resulted from the inclusion 
of the potential factor $\frac {1}{\Om_x}$, is consistent with 
this symmetry.

The 3D Euler equation can be written in vorticity form,
\begin{equation}
\pa_t \Om + (u \cdot \na) \Om - (\Om \cdot \na) u = 0 \ ,
\label{3deuler}
\end{equation}
where $u = (u_1, u_2, u_3)$ is the velocity, $\Om = (\Om_1, \Om_2, \Om_3)$
is the vorticity, $\na = (\pa_x, \pa_y, \pa_z)$, 
$\Om = \na \times u$, and $\na \cdot u = 0$. $u$ can be 
represented by $\Om$ for example through Biot-Savart law.

{\bf Theorem 7.} \cite{LY03} \it The Lax pair of the 3D Euler equation 
(\ref{3deuler}) is given as
\begin{equation}
\left \{ \begin{array}{l} 
L \phi = \la \phi \ ,
\\
\pa_t \phi + A \phi = 0 \ ,
\end{array} \right.
\label{alaxpair}
\end{equation}
where
\[
L \phi = (\Om \cdot \na )\phi \ , 
\ \ \ A \varphi = (u \cdot \na )\phi \ , 
\]
$\la$ is a complex constant, and $\phi$ is a complex scalar-valued function. \rm

Our hope is that the infinitely many conservation laws generated by $\la 
\in C$ can provide a priori estimates for the global well-posedness of 
3D Navier-Stokes equations, or better understanding on the global 
well-posedness. 

Notice that the equation (\ref{3deuler}) without the constraint 
$\Om = \na \times u$ and $\na \cdot u = 0$ is a compatibility condition of 
the Lax pair (\ref{alaxpair}). In general, Lax pairs can support finite time blow up 
solutions. Our conjecture here is that adding only the first constraint $\Om = \na \times u$,
the equation (\ref{3deuler}) may have finite time blow up solutions which will be removed by 
the second constraint $\na \cdot u = 0$. 

\subsection{Invariant Manifolds and Their Zero-Viscosity Limits}

To begin a dynamical system study on Navier-Stokes equations, one needs to study the spectra of 
their linearizations at fixed points and their invariant manifolds. 

Following the notations of (\ref{euler}) and (\ref{3deuler}), we shall study the Navier-Stokes 
equations in the forms:
\begin{equation}
\pa_t \Om + \{ \Psi, \Om \} = \nu [\Dl \Om + f(x)] \ ,
\label{2DNS}
\end{equation}
for 2D, and 
\begin{equation}
\pa_t \Om + (u \cdot \na) \Om - (\Om \cdot \na) u = \nu [\Dl \Om + f(x)] \ ,
\label{3DNS}
\end{equation}
for 3D where $\nu$ is the viscosity, $\Dl$ is the Laplacian, and $f(x)$ is 
the body force. Periodic boundary condition is imposed.

{\bf Theorem 8.} \cite{Li05} \it When $\nu > 0$, any fixed point of NS has invariant 
manifolds. The provable size of the invariant manifolds shrinks to zero as $\nu \ra 0^+$. \rm

The interesting open question is: Can the size of some invariant manifold be $\cO (1)$ 
as $\nu \ra 0^+$ ?

For simple examples of fixed points, detailed information on the spectra of the corresponding 
linear NS operator can be obtained \cite{Li05} \cite{Li00} \cite{LLS04}. In the 2D case, consider 
the simple shear $\Om_* = \cos y$ defined on the rectangular periodic domain 
$[0, 2\pi /\al ] \times [0, 2\pi ]$ where $1/2 < \al < 1$. The corresponding 
linear NS operator decouples into infinitely many sub-operators labeled by ($\hk_1 , \hk_2$) $\in 
\ZZo$.

{\bf Theorem 9.} \cite{Li05} \it The spectra of the 2D linear NS operator have the following 
properties labeled by ($\hk_1 , \hk_2$):
\begin{enumerate}
\item $(\al \hk_1)^2+(\hk_2+n)^2 > 1$, $\forall n \in \Z$. 
When $\nu \ra 0$, there is no eigenvalue of non-negative real part. 
When $\nu = 0$, the entire spectrum is the continuous spectrum
\[
\left [ -i\frac{\al |\hk_1|}{2}, \ i\frac{\al |\hk_1|}{2} \right ]\ .
\]
\item $\hk_1 = 0$, $\hk_2 = 1$. The spectrum consists of the eigenvalues 
\[
\la = - \nu n^2 \ , \quad n \in \Zo \ .
\]
The eigenfunctions are the Fourier modes
\[
\tom_{np} e^{inx_2} + \ \mbox{c.c.}\ \ , \quad \forall \tom_{np} \in 
\C\ , \quad n \in \Zo \ .
\]
As $\nu \ra 0^+$, the eigenvalues are dense on the negative half of the real 
axis $(-\infty, 0]$. Setting $\nu =0$, the only eigenvalue is $\la = 0$ of 
infinite multiplicity with the same eigenfunctions as above.
\item $\hk_1 = -1$, $\hk_2 = 0$. (a). $\nu >0$. For any $\al \in (0.5, 0.95)$,
there is a unique $\nu_*(\al)$,
\[
\frac{\sqrt{32-3\al^6-17\al^4-16\al^2}}{4(\al^2+1)(\al^2+4)} < 
\nu_*(\al) < \frac{1}{2(\al^2+1)} \sqrt{\frac{1-\al^2}{2}}\ ,
\]
where the term under the square root on the left is positive for 
$\al \in (0.5, 0.95)$, and the left term is always less than the right term.
When $\nu > \nu_*(\al)$, there is no eigenvalue of non-negative real part. 
When $\nu = \nu_*(\al)$, $\la =0$ is an eigenvalue, and all the rest 
eigenvalues have negative real parts. When $\nu < \nu_*(\al)$, there is 
a unique positive eigenvalue $\la (\nu )>0$, and all the rest 
eigenvalues have negative real parts. $\nu^{-1} \la (\nu )$ is a strictly 
monotonically decreasing function of $\nu$. When $\al \in (0.5, 0.8469)$,
we have the estimate
\begin{eqnarray*}
& & \sqrt{\frac{\al^2(1-\al^2)}{8(\al^2+1)}-\frac{\al^4 (\al^2+3)}{16
(\al^2+1)(\al^2+4)}} - \nu (\al^2+1) < \la (\nu ) \\
& & < \sqrt{\frac{\al^2(1-\al^2)}{8(\al^2+1)}}- \nu \al^2 \ ,
\end{eqnarray*}
where the term under the square root on the left is positive for 
$\al \in (0.5, 0.8469)$.
\[
\sqrt{\frac{\al^2(1-\al^2)}{8(\al^2+1)}-\frac{\al^4 (\al^2+3)}{16
(\al^2+1)(\al^2+4)}} \leq \lim_{\nu \ra 0^+} \la (\nu )  \leq 
\sqrt{\frac{\al^2(1-\al^2)}{8(\al^2+1)}} \ .
\]
In particular, as $\nu \ra 0^+$, $\la (\nu ) =O(1)$.

(b). $\nu =0$. When $\al \in (0.5, 0.8469)$, we have only two eigenvalues
$\la_0$ and $-\la_0$, where $\la_0$ is positive,
\[
\sqrt{\frac{\al^2(1-\al^2)}{8(\al^2+1)}-\frac{\al^4 (\al^2+3)}{16
(\al^2+1)(\al^2+4)}} < \la_0 <
\sqrt{\frac{\al^2(1-\al^2)}{8(\al^2+1)}} \ .
\]
The rest of the spectrum is a continuous spectrum $[-i\al /2, \ i\al /2]$.

(c). For any fixed $\al \in (0.5, 0.8469)$,
\[
\lim_{\nu \ra 0^+} \la (\nu ) = \la_0 \ .
\]
\item Finally, when $\nu = 0$, the union of all the above pieces of 
continuous spectra is the imaginary axis $i\R$.
\end{enumerate} \rm

In general, for any fixed point of NS, the spectrum of the linear NS operator 
consists of only eigenvalues which lie in a parabolic region \cite{Li05}. The 
above theorem establishes some properties of the zero-viscosity limits of these 
eigenvalues for a simple shear. From numerical simulations on many examples \cite{LL07}, 
these eigenvalues often undergo fascinating deformations in the zero-viscosity limits,
which can be classified into the following four categories:
\begin{enumerate}
\item {\em Persistence:} These are the eigenvalues that persist and approach 
to the eigenvalues of the corresponding linear Euler operator when the viscosity 
approaches zero. (e.g. at 2D and 3D shears, and cat's eye.)
\item {\em Condensation:} These are the eigenvalues that approach and form 
a continuous spectrum for the corresponding linear Euler operator when the viscosity 
approaches zero. (e.g. at 2D and 3D shears, cat's eye, and ABC flow.)
\item {\em Singularity:} These are the eigenvalues that approach to a 
set that is not in the spectrum of the corresponding linear Euler operator when 
the viscosity approaches zero. (e.g. at 2D and 3D shears.)
\item {\em Addition:} This is a subset of the spectrum of the linear Euler operator,
which has no overlap with the zero viscosity limit set of the spectrum of the linear 
NS operator. (e.g. cat's eye.)
\end{enumerate} 

Investigating the zero-viscosity limit is very important for studying the NS dynamics 
near the Euler regime. By a combination of numerics and analysis in \cite{LL07}, 
we have tried to push the Melnikov integral tool into the study on chaotic dynamics of 
NS. One interesting question is the conjecture that Euler equations have heteroclinics \cite{LL07}.

Proving the existence of chaos in NS will provide a foundation for control of turbulence.
Control of turbulence depends on the information that supports the existence of chaos. 
In fact, control can be regarded as an extra boby force in the NS equations. Thus, 
control of chaos (turbulence) is a sub-problem of existence of chaos (turbulence). 
From this perspective, investigating the existence of chaos in NS has great industrial value.

\subsection{Global Regularity of Navier-Stokes Equations}

A notable mathematical problem on 3D Navier-Stokes equations is the global regularity problem
which was chosen as one of the seven Clay problems. Specifically, the fact \cite{Ler34} that 
\[
\int \int | \na u|^2 \ dx\ dt  
\]
being bounded only implies 
\[
\int | \na u|^2 \ dx
\]
being bounded for almost all $t$, is the key of the difficulty. In fact, 
Leray was able to show that the possible exceptional set of $t$ is 
actually a compact set of measure zero. There have been a lot of more 
recent works on describing this exceptional compact set \cite{CKN82}. 
The claim that this possible exceptional compact set is actually empty, 
will imply the global regularity and the solution of the problem. The 
hope for such a claim seems slim. There are weaker conditions implying 
global regularity, e.g. $L^4$ norm of velocity being bounded for t in the exceptional 
set \cite{Con01} \cite{Yud03}. But I think the weak and the strong conditions will be true 
or false simultaneously. There are several alternative ways 
of approaching the global regularity problem \cite{Tao07} including 
the Lax pair for 3D Euler equations in a previous subsection.

Even for ordinary differential equations, often one can not prove their 
global well-posedness, but their solutions on computers look perfectly 
globally regular and sometimes chaotic. Chaos and global regularity are 
compatible. The hallmark of chaotic solutions is their sensitive dependence 
on initial conditions. That is, a small change in the initial conditions 
will lead to a huge change after sufficiently long time, while all the 
solutions are kept in a bounded region in the phase space. This leads to the 
unpredictability, not irregularity, of chaotic solutions.
The fact that fluid experimentalists quickly discovered 
shocks in compressible fluids and never found any finite time blow up in 
incompressible fluids, indicates that there might be no finite time 
blow up in 3D Navier-Stokes equations (even Euler equations). On 
the other hand, the solutions of 3D Navier-Stokes equations can definitely be 
turbulent. 

Replacing the viscous term $\nu \ \Dl u$ by higher order 
derivatives, one can prove the global regularity \cite{KP02}.
This leaves the global regularity of a more challenging and 
interesting mathematical problem. Assume that the unthinkable event 
happens, that is, someone proves the existence of a meaningful finite 
time blow up in 3D NS, then fluid experimentalists need to identify 
such a finite time blow up in the experiments. If they fail, then the choice 
will be whether or not to replace the viscous term 
$\nu \ \Dl u$ in the Navier-Stokes equations
by higher order derivatives to better model the fluid motion.
Such a freedom was allowed in the original derivation of NS. Specifically, 
the fluid shear stress can depend on second or higher derivatives of velocity
besides the fluid strain given by first derivatives of velocity. 

The global regularity and turbulence are separate problems. Turbulence is more 
of a dynamical system problem. Of course, global attractors depend on global 
well-posedness \cite{CV02}. Often a detailed local dynamical system 
study does not depend on 
global well-posedness. Local well-posedness is often enough. In fact, this 
is the case in my proof on the existence of chaos in partial differential 
equations \cite{Li04}.

\subsection{Searching for an Effective Description of Chaotic (Turbulent) Solutions Beyond 
the Reynolds Average}

Of the greatest industrial value is an effective description of chaotic (turbulent) solutions 
beyond the Reynolds average. From what we learn about chaos in partial 
differential equations \cite{Li04}, turbulent 
solutions not only have sensitive dependences on initial conditions, but 
also are densely packed inside a domain in the phase space. They are far 
away from the feature of fluctuations around a mean. In fact, they wander 
around in a fat domain rather than a thin domain in the phase space. 
Therefore, averaging makes no sense at all. One has to seek other descriptions. 
In \cite{Li07a}, an interesting ``segment description'' is proposed.

\end{document}